\documentclass[12pt,letterpaper]{amsart} 
\usepackage{epic,eepic,latexsym, amssymb, amscd, amsfonts, xypic, graphicx}

 \newlength{\baseunit}               
 \newcount{\numlines}                
\newcommand{\tpoint}[1]{\vspace{3mm}\par \noindent \refstepcounter{subsection}{\bf \thesubsection.}
  {\em #1. ---} }
\newcommand{\epoint}[1]{\vspace{3mm}\par \noindent \refstepcounter{subsection}{\bf \thesubsection.}
  {\em #1.} }
\newcommand{\bpoint}[1]{\vspace{3mm}\par \noindent \refstepcounter{subsection}{\bf \thesubsection.}
  {\bf #1.} }

\newtheorem{tm}{Theorem}[subsection]

\newtheorem{lm}[tm]{Lemma}

\newcommand{\Z}{\mathbb{Z}}

\newcommand{\C}{\mathbb{C}}

\newcommand{\proj}{\mathbb P}

\newcommand{\al}{\alpha}

\newcommand{\De}{\Delta}
\newcommand{\si}{\sigma}

\newcommand{\Pic}{\operatorname{Pic}}

\newcommand{\Div}{\operatorname{Div}}
\newcommand{\Aut}{\operatorname{Aut}}
\newcommand{\Def}{\operatorname{Def}}
\newcommand{\Gal}{\operatorname{Gal}}

\newcommand{\secretnote}[1]{}

\newcommand{\lremind}[1]{{}}

\newcommand{\Q}{\mathbb{Q}}
\newcommand{\Qbar}{\overline{\mathbb{Q}}}

\newcommand{\Hom}{\operatorname{Hom}}

\newcommand{\comment}[1]{}

\begin{document}
\pagestyle{plain}

\title{\large{
Absolute
Galois acts faithfully on the components of the moduli
space of surfaces:  A Belyi-type theorem in higher dimension}}

\date{April 12, 2007.}

\author{Robert W. Easton and Ravi Vakil}

\address{Department of Mathematics, Stanford University, Stanford, California 94305}
\email{easton@math.stanford.edu \\ vakil@math.stanford.edu}

\begin{abstract}
  Given an object over $\Qbar$, there is often no reason for
  invariants of the corresponding holomorphic object to be preserved
  by the absolute Galois group $\Gal(\Qbar/\Q)$, and in general this
  is not true, although it is sometimes surprising to observe in
  practice.  The case of covers of the projective line branched only
  over the points $0$, $1$, and $\infty$, through Belyi's theorem,
  leads to Grothendieck's {\em dessins d'enfants} program for
  understanding the absolute Galois group through its faithful action
  on such covers.  This note is motivated by Catanese's question about
  a higher-dimensional analogue: does the absolute Galois group act
  faithfully on the deformation equivalence classes of smooth
  surfaces? (These equivalence classes are of course by definition the
  strongest deformation invariants.)  We give a short proof of a
  weaker result: the absolute Galois group acts faithfully on the {\em
    irreducible} components of the moduli space of smooth surfaces (of
  general type, canonically polarized).   Bauer, Catanese,
  and Grunewald have recently answered Catanese's original question
  using a different construction \cite{bcgnew}.
\end{abstract}
\maketitle


\vspace{-0.2in}
{\parskip=12pt 

\section{Introduction}

Given a object defined over $\Qbar$, certain topological invariants of
the corresponding holomorphic object are known to be preserved by the
absolute Galois group $\Gal(\Qbar/\Q)$. This is because these
invariants are algebraic in nature.  For example, if $X$ is a
nonsingular projective variety, the Betti numbers are algebraic (shown
by Serre in his GAGA paper, \cite{serre1}). The profinite completion
of the fundamental group of $X$ is the \'etale fundamental group.
More generally, Artin and Mazur showed that the profinite completion
of the homotopy type of $X$ is algebraic \cite{am}.

It is thus natural to ask what topological invariants of
the corresponding holomorphic object are preserved by conjugation.
Indeed, given an object defined over $\Qbar$, there is often no reason
for topological invariants of the corresponding holomorphic object to
be preserved by the absolute Galois group.  In the
case of covers of the projective line branched only over the points
$0$, $1$, and $\infty$, this leads to Grothendieck's {\em dessins
  d'enfants} program for understanding the absolute Galois group
\cite{Grothendieck}, through its faithful action on such covers.  In
other words, given any nontrivial element $\si \in \Gal(\Qbar/\Q)$, there is a
cover $C \rightarrow \proj^1$ (over $\Qbar$) such that $\si(C)
\rightarrow \proj^1$ is a topologically different cover (where both
covers are now considered over $\C$, as maps of Riemann surfaces).

Similarly, returning to the case of smooth varieties, Serre gave an
elegant example \cite{S} of a smooth variety $X$ over $\Qbar$ and an
element $\si \in \Gal(\Qbar/\Q)$ such that the fundamental groups of
the complex manifolds $X$ and $\si(X)$ are different.  (As the
profinite completions $\pi_1^{\text{\'et}}(X)$ and
$\pi_1^{\text{\'et}}(\si(X))$ are isomorphic, the fundamental groups
are necessarily infinite.)  Abelson \cite{A} gave examples of
conjugate (nonsingular projective) varieties with the {\em same}
fundamental group yet of different homotopy types.  He also gave
examples of conjugate (nonsingular quasiprojective) varieties that are
homotopy equivalent but not homeomorphic.  More examples of
nonhomeomorphic conjugate varieties have been given quite recently by
Artal Bartolo, Carmona Ruber, and Cogolludo Agust\'{i}n \cite{three},
and Shimada \cite{shimada}.

Surprising examples of a different flavor, using Beauville surfaces,
were given by Catanese earlier (see Theorem 21 and the discussion
just before Question 4 in \cite{BCP}; cf.\ \cite[Thm.~3.3]{Cat03} and
\cite[Thm.~4.14]{Cat00}). 

A potentially rich third family of examples arises from the theory of
Shimura varieties, as described by Milne \cite[p.~7]{milne}.  By a
theorem of Baily and Borel \cite{bb}, the quotient of a bounded
symmetric domain by an arithmetic subgroup of its analytic
automorphism group has a canonical structure of a quasiprojective
complex variety $V$.  A conjecture of Langlands implies that if $\si
\in \Gal(\Qbar/\Q)$, then $\si(V)$ is again such a quotient, and
describes explicitly what the bounded symmetric domain and arithmetic
subgroup are; this conjecture was proved by Borovoi and Milne
\cite{borovoi, m1} using a theorem of 
Kazhdan and Nori-Raghunathan  \cite{k1, k2, nr}.  One
should be able to show that these arithmetic groups (the fundamental
groups of the Shimura varieties in cases of good quotients) are not
isomorphic (as abstract groups); to our knowledge, the details have
not yet been worked out in the literature.

The strongest deformation-invariant discrete invariant is of course
the deformation equivalence class.  This note is motivated by a
question of Catanese: does the absolute Galois group act faithfully on
the deformation equivalence classes of surfaces (defined over
$\Qbar$)?  In other words, given any nontrivial $\si \in
\Gal(\Qbar/\Q)$, can one produce a surface $X$ such that $\si(X)$ is
not deformation-equivalent to $X$?  Catanese has shown that it is
indeed true when $\si$ is complex conjugation
(see \cite[Thm.~3.5]{Cat03}, as well as later numerous rigid examples by
Bauer, Catanese, and Grunewald in \cite{bcg2}).  (One might speculate that
every element of $\Gal(\Qbar/\Q)$ other than the identity and complex
conjugation can change the homeomorphism type of a $\Qbar$-variety,
and combined with Catanese's example this would answer Catanese's
question.  However, it is not clear from the examples produced to date
that {\em every} $\si \in \Gal(\Qbar/\Q)$ besides identity and complex conjugation has this property.)

 In effect, Catanese's question translates to:

 \epoint{Catanese's question} {\em Does $\Gal(\Qbar/\Q)$ act
   faithfully on the connected components of the moduli space of
   surfaces of general type?}
\label{question}

We are able to show the following weaker result:

\tpoint{Main Theorem}  {\em  The absolute Galois group acts
faithfully on the {\bf irreducible} components of the moduli space
of surfaces of general type. More precisely, for {\bf any} nontrivial $\si
\in \Gal(\Qbar/\Q)$, we exhibit a surface $X$ over $\Qbar$, where
$X$ has ample canonical bundle (indeed very ample), and such that
$X$ and $\si(X)$ do not lie on the same component of the moduli
space.}

\noindent {\em Important remark.}  After we wrote this note, Bauer,
Catanese, and Grunewald informed us that they have given a complete
answer to Catanese's question using a different construction, and
outlined their ingenious argument. Catanese announced their work at
the Alghero conference of September 2006, and the paper will be public
very shortly \cite{bcgnew}.

\noindent {\bf Strategy.}  We first choose $z \in \Qbar$ not fixed by
$\si$.  Our surface $X=X_z$ will be constructed so that the number $z$
will be ``encoded'' in it (and its infinitesimal deformations), and
such that its conjugate $\si(X) = X_{\si(z)}$ (and its infinitesimal
deformations) will encode the number $\si(z)$ in the same way.  Thus
there are {\em Zariski} neighborhoods of the points $[X_z]$ and
$[X_{\si(z)}]$ of the moduli space that are disjoint.

We perform this encoding by first describing a configuration of points
and lines on the plane (over $\Qbar$) such that the combinatorics of
incidences of points and lines encodes the number $z$, in such a way
that the $\si$-conjugate encodes the number $\si(z)$.  We do this as
follows. There will be four distinguished ordered points on a line in
the plane; they will be the four points with the most lines through
them. The cross-ratio of these four points on the line will be $z$.
Hence the Galois conjugate will have cross-ratio $\si(z)$.  

Then (as in \cite{V}) we let $X$ be a branched cover of the blow-up
of the plane at the points, where the branch locus consists of the
proper transforms of the lines, as well as several high-degree
curves. This positivity will force the vanishing of certain
cohomology groups, which will allow us to ensure that the
deformations of $X$ correspond exactly with the deformations of the
point-line configuration on the plane.  More precisely, from $X$ (or
any infinitesimal deformation), we can recover the branched cover,
and hence the data of the point-line configuration.  These
constructions ``commute with $\si$'', yielding the result.

\bpoint{\bf Miscellaneous remarks} {\em (a)} We do not know if the
two surfaces are homeomorphic (if $\sigma$ is not complex
conjugation), and we have no reason to expect that they are. If they are
not homeomorphic, then this would answer Catanese's question ~\ref{question}
in the
affirmative.  Moreover, they are constructed so that it is possible in
theory to compute their fundamental groups.  If one could do so, and
show that they are different, this would answer Catanese's question
completely.  

{\em (b)} Gonz\'alez-Diez \cite{G} and Paranjape \cite{paranjape}
have given other higher-dimensional analogues of Belyi's theorem.

{\em (c)} This is vaguely reminiscent of dessins d'enfants.  In the
case of covers of $\proj^1$, the covers were encoded by graphs --- not
just the incidences of vertices and edges, but also the embedding in
the surface.  In this case, the surfaces are encoded by lines in the
plane --- not just the combinatorial data of incidences of points and
lines, but also the embedding in the (complex) plane.

\noindent {\bf Acknowledgments.}
This article was motivated by a lecture of Ingrid Bauer at
Oberwolfach, and a seminar series by Greg Brumfiel on dessins
d'enfants at Stanford, and we are grateful to them both for
introducing us to these questions.  We also thank Matt Emerton, Soren
Galatius, James Milne, Mihran Papikian, Kapil Paranjape, and Dinakar
Ramakrishnan for suggestions and references.  We are grateful to
Bauer, Fabrizio Catanese, and Fritz Grunewald for informing us about
\cite{bcgnew} and explaining their argument.

\section{The argument}

Take any nontrivial $\sigma \in \Gal(\bar{\Q}/\Q)$, and fix some $z
\in \bar{\Q}$ with $\sigma(z)\neq z$.\label{theargument}

\bpoint{Point-line configurations} We construct a point-line
configuration encoding the algebraic number $z$ as follows.  Let
$p(z)$ be the minimal polynomial of $z$ (with relatively prime integer
coefficients).  Choose a line $\ell$ on the plane, and three distinct
points on it, which we name $0$, $1$, $\infty$. Then the points of
$\ell - \{\infty\}$ are naturally identified with numbers; i.e.,
elements of $\bar{\Q}$. Given three points $a, b$ and $c$ on
$\ell-\{\infty\}$ (considered as numbers), it is straightforward to
construct a point-line configuration through these points that forces
precisely the equation $a+b=c$; and a different configuration that
forces $ab=c$; and a third that forces $a=-b$. We combine these
operations suitably so as to force $p(z)=0$.  (This sort of recipe is
well-known and straightforward, so we omit the details.  See for
example \cite[p.~13]{shafarevich}.  The apotheosis of this idea is
Mn\"ev's Universality Theorem \cite{M1, M2}.)

Let $\mathcal{L}_z'$ denote this configuration representing $z$. We
modify $\mathcal{L}_z'$ so as to produce a configuration
$\mathcal{L}_z \supset \mathcal{L}_z'$ over which a branched
$(\mathbb{Z}/2)^3$-cover is readily constructed.  First, we add a
general line through each point in the plane through which an odd
number of lines (greater than one) in $\mathcal{L}_z$ pass.
Finally, we add general lines through the points $0$, $1$, $\infty$,
and $z$ (an even number through each), so that the points on the
plane with the most lines through them are, in order, these four
points. The points in our point-line configuration will consist of
all points of intersection of pairs of such lines.  Note that from a
general such configuration, we can recover $z$ by finding the four
points with the most lines through them, observing that they lie on a
common line, and taking their cross-ratio on this line.  Also,
acting on such a configuration with $\si$ will yield a configuration
encoding $\si(z)$ in the same way. Thus the first point-line
configuration may not be deformed to the second (while preserving the
point-line incidences). The ``even valence'' condition will be used
later.

\bpoint{Branched covers background}
This portion of our
argument follows \cite{V}.
We will obtain $X$ by blowing up the plane at our marked points, and
taking a suitable branched cover of the resulting surface.

We first review some results about branched covers, due to Catanese,
Pardini, Fantechi, and Manetti. Suppose $G=(\Z/p\Z)^n$ with $p$
prime.  Let $G^{\vee}= \Hom(G,\mathbb{C}^*)$ be the group of complex
characters of $G$, and for each $\chi \in G^{\vee}$, define
$(\chi,g)\in \{0,\ldots,p-1\}$ by
$\chi(g)=e^{\frac{2\pi i}{p}(\chi,g)}$. Let $S$ be any nonsingular
surface, and suppose $\{D_g\}_{g \in G}, \{M_{\chi}\}_{\chi
\in G^{\vee}}$ are divisors in $S$ satisfying $D_0=\emptyset$ and
\[
pM_{\chi}\equiv \underset{g \in G}{\sum} (\chi, g)D_g
\]
in $\Pic(S)$ for all $\chi \in G^{\vee}$.  Moreover, suppose the
$D_g$ are all nonsingular curves, no three intersect in a point, and
$D_g\cap D_{g'} \neq \emptyset$ only if $g$ and $g'$ are independent
in $G$ (i.e.\ $g'$ is not a multiple of $g$). Then:

\tpoint{Theorem}{\em  \label{PropositionA} There exists a
nonsingular $G$-cover $\pi \colon X\to S$ with branch divisor
$D=\underset{g \in G}{\bigcup}D_g$. Moreover, if $n\geq 3$ and
$M_{\chi}$ is sufficiently ample for all nonzero $\chi \in
G^{\vee}$, then:
\begin{enumerate}
\item[(a)] $K_X$ is very ample;
\item[(b)] deformations of $(S, \{D_g\})$ are equivalent to
  deformations of $X$, i.e.\ the natural map $$\Def (S, \{D_g\})
  \rightarrow \Def(X)$$ is an isomorphism; and
\item[(c)] $\Aut(X)\cong G$.
\end{enumerate}
}

Part (a) is given in the e-print version of \cite{V} (Theorem 4.4),
and the idea is due to Catanese.  (The argument for bidouble covers is
given in \cite[p. ~502]{C}.)  Part (b) is \cite[Thm.\ 4.4(c)]{V}, and
the argument is due to Manetti \cite[Cor.\ 3.23]{manetti}; indeed, the
case $p=2$ that we will actually use is Manetti's original result.
Part (c) is due to Fantechi and Pardini \cite[Thm.\ 4.6]{FP}.

\bpoint{Branched covers of the blow-up of the plane} We now let $S_z$
be the blow-up of $\proj^2$ at our marked points, and $C_z$ be the
strict transform of the union of our lines. We will construct a
branched cover with $G = (\Z / 2 \Z)^3$.  First, we define maps
$D:G\to \Div(S_z)$ and $M:G^{\vee}\to \Pic(S_z)$. Let $D_0=\emptyset$.
Fix any nonzero $\al\in G$ and let $D_{\al}=C_z$.  Let $L$ denote the
number of lines in our configuration.  Fix any map $m:G\to
\mathbb{Z}^+$ such that $m_0=0, m_{\al}=L$, and $\sum_{g\in G}m_gg=0$
in $G$.  Then $\sum_{g\in G}(\chi, g)m_g$ is even for every $\chi \in
G^{\vee}$:  
$$1= \chi(0) = \chi \left( \sum_g m_g g \right) = \prod_g \chi(g)^{m_g} = 
e^{\frac {2 \pi i} 2 \sum_g \left(
\chi, g \right) m_g}
= (-1)^{ \sum_g \left(
\chi, g \right) m_g}.$$ 
We will use this fact momentarily.

For $g\in G-\{0,\al\}$, define $D_g$ to be the pull-back of a
general (nonsingular) curve in $\proj^2$ of degree $m_g$. By our choice of
$m_g$, we then have
\begin{align*}
\underset{g \in G}{\sum}(\chi, g)D_g &\equiv
-(\chi,\al)\underset{q}{\sum}e_q(\mathcal{L}_z)E_q+\underset{g\in
G}{\sum}(\chi,g)m_gH,
\end{align*}
in $\Pic(S_z)$, where $H$ is the hyperplane class in $\proj^2$ and $e_q(\mathcal{L}_z)$ is the number of lines
in our configuration passing through the point $q$. By construction,
the above divisor is even, and hence we may define $\chi:
G^{\vee}\to \Pic(S_z)$ by $M_{\chi}=\frac{1}{2}\sum_g (\chi, g)D_g$.
Note that the $M_{\chi}$ can be made arbitrarily ample by an
appropriate choice of the map $m$. For such a choice, by
Theorem~\ref{PropositionA} we obtain a nonsingular general type
$G$-cover $\pi\colon X_z\to S_z$ with branch divisor $D=\bigcup
D_g\supset C_z$. The same construction mutatis mutandis produces a
conjugate nonsingular general type $G$-cover $\pi_{\sigma(z)}\colon
X_{\sigma(z)}\to S_{\sigma(z)}$ with branch divisor
$D_{\sigma(z)}=\sigma(D_z)\supset \sigma(C_z)=C_{\sigma(z)}$. 
It then follows from
Theorem~\ref{PropositionA}(b) that the deformations of $X_z$ (resp.\
$X_{\si(z)}$) are equivalent to the deformations of $(S_z,D_z)$ (resp.\
$(S_{\sigma(z)},D_{\sigma(z)})$).

We now describe how to recover the number $z$ from $S_z$ and any
infinitesimal deformation (and similarly for $S_{\si(z)}$).  By
Theorem~\ref{PropositionA}(c), $G \rightarrow \Aut(X_z)$ is an
isomorphism, from which we may recover $X_z \rightarrow X_z / G =
S_z$.  The components of the branch divisor of $X_z \rightarrow X_z/
G$ are the divisors $\{D_g\}_{g\neq 0}$.  All but one of them
(all except $D_{\al}$) are $\Q$-multiples of each other; they are all equivalent
to multiples of $H$, the pullback of the hyperplane divisor in
$\proj^2$.  We may therefore use any divisor from this distinguished
collection to recover the blow-down to $\proj^2$.  The components of
the remaining branch divisor $D_{\al}$ give the lines
$\mathcal{L}_z\subset \proj^2$.  From this we may recover $z$.  This
discussion clearly extends to the (formal) deformation space around
$[X_z]$ in the moduli space.  (The argument begins: Let $\De$ be the
formal deformation space, and $\mathcal{X}_z$ be the total family of
the deformation.  Then $\Aut( \mathcal{X}_z / \De)$ is the trivial
group scheme $G$ over $\De$, from which we obtain $\mathcal{X}_z
\rightarrow \mathcal{S}_z$, etc.)


\bpoint{Closing remarks}  This result suggests an approach to answering
Catanese's question~\ref{question} in general, by producing a rigid
surface as such a branched cover.  One might attempt to do so by
rigidifying the point-line configuration $\mathcal{L}'_z$ of the start
of \S \ref{theargument} by adding judiciously chosen additional lines,
using a theorem of Paranjape \cite[Thm.~2]{paranjape}.  One would then
have to modify the argument to ensure that {\em (a)} the four points
``marking $z$'' remain distinguished, and {\em (b)}
Theorem~\ref{PropositionA}  continues to hold, without the assistance of
the positivity of $M_{\chi}$.

} 

\end{document}